\documentclass[12 pt]{article}
\usepackage{amsmath,amssymb,amscd,amsthm,verbatim}
\input epsf                 

\newtheorem{proposition}{Proposition}[section]
\newtheorem{theorem}[proposition]{Theorem}
\newtheorem{corollary}[proposition]{Corollary}

\newtheorem{prop}[proposition]{Proposition}

\newtheorem{cor}[proposition]{Corollary}

\newtheorem{conj}[proposition]{Conjecture}

\newcommand{\naturals}{\mathbb N}

\newcommand{\integers}{\mathbb Z}

\newcommand{\reals}{\mathbb R}

\newcommand{\F}{{\mathcal F}}
\newcommand{\Po}{{\mathcal P}}

\newcommand{\covers}{{\,\cdot\!\!\!\! >}}

\newcommand{\rank}{{\mathrm{rank}}}

\newcommand{\0}{{\hat{0}}}
\newcommand{\1}{{\hat{1}}}
\newcommand{\set}[1]{{\lbrace #1 \rbrace}}

\newcommand{\join}{\vee}

\newcommand{\meet}{\wedge}

\newcommand{\Pyr}{{{\mathrm{Pyr}}}}
\newcommand{\Sh}{{{\mathrm{Sh}}}}
\newcommand{\Z}{{{\mathrm{Zip}}}}

\title{The cd-index of Bruhat intervals}

\author{Nathan Reading
\thanks{The author was partially supported by the Thomas H. Shevlin Fellowship from the University of Minnesota Graduate School and by 
NSF grant DMS-9877047.  This article contains material from the author's doctoral thesis~\cite{thesis}.}\\
\small Mathematics Department\\[-0.8ex]
\small University of Michigan, Ann Arbor, MI 48109-1109\\[-0.8ex]
\small \texttt{nreading@umich.edu}}

\date{\small MR Subject Classifications: 20F55, 06A07}

\begin{document}
\maketitle

\begin{abstract}
We study flag enumeration in intervals in the Bruhat order on a Coxeter group by means of a structural recursion on intervals in the Bruhat order.
The recursion gives the isomorphism type of a Bruhat interval in terms of smaller intervals, using basic geometric operations which preserve 
PL sphericity and have a simple effect on the cd-index. 
This leads to a new proof that Bruhat intervals are PL spheres as well a recursive formula for the cd-index of a Bruhat interval.
This recursive formula is used to prove that the cd-indices of Bruhat intervals span the space of cd-polynomials.

The structural recursion leads to a conjecture that Bruhat spheres are ``smaller'' than polytopes.
More precisely, we conjecture that if one fixes the lengths of $x$ and $y$, then the cd-index of a certain dual stacked polytope is a coefficientwise 
upper bound on the cd-indices of Bruhat intervals $[x,y]$.
We show that this upper bound would be tight by constructing Bruhat intervals which are the face lattices of these dual stacked polytopes.
As a weakening of a special case of the conjecture, we show that the flag h-vectors of lower Bruhat intervals are bounded above by the flag 
h-vectors of Boolean algebras (i.\ e.\ simplices).
\end{abstract}

A graded poset is {\em Eulerian} if in every non-trivial interval, the number of elements of odd rank equals the number of elements of even rank. 
Face lattices of convex polytopes are in particular Eulerian and the study of flag enumeration in Eulerian posets has its origins in the 
face-enumeration problem for polytopes.
All flag-enumerative information in an Eulerian poset $P$ can be encapsulated in a non-commutative generating function $\Phi_P$ 
called the {\em cd-index}.
The cd-indices of polytopes have received much attention, for example in \cite{Bayer-Billera,Bayer-Klapper,Ehr-Fox,Kalai,Stanley2}.

A {\em Coxeter group} is a group generated by involutions, subject to certain relations.
Important examples include finite reflection groups and Weyl groups.
The {\em Bruhat order} on a Coxeter group is a partial order which has important connections to the combinatorics and representation theory 
of Coxeter groups, and by extension Lie algebras and groups. 
Intervals in Bruhat order comprise the other important class of Eulerian posets.
However, flag enumeration for Bruhat intervals has previously received little attention.
The goal of the present work is to initiate the study of the cd-index of Bruhat intervals.

The basic tool in our study is a fundamental structural recursion (Theorem~\ref{etathetazippers}) on intervals in the Bruhat order on 
Coxeter groups.
This recursion, although developed independently, has some resemblance to work by du Cloux \cite{DuCloux} and by Dyer \cite{Dyerthesis}.
The recursion gives the isomorphism type of a Bruhat interval in terms of smaller intervals, using some basic geometric operations, 
namely the operations of pyramid, vertex shaving and a ``zipping'' operation.
The result is a new inductive proof of the fact~\cite{Bj} that Bruhat intervals are PL spheres (Corollary~\ref{spheres}) as well as recursions for the 
cd-index of Bruhat intervals (Theorem~\ref{formula}).

The recursive formulas lead to a proof that the cd-indices of Bruhat intervals span the space of cd-polynomials (Theorem~\ref{span}), 
and motivate a conjecture on the  upper bound for the cd-indices of Bruhat intervals (Conjecture~\ref{dual stacked upper}).
Let $[u,v]$ be an interval in the Bruhat order such that the rank of $u$ is $k$ and the rank of $v$ is $d+k+1$.
We conjecture that the coefficients of $\Phi_{[u,v]}$ are bounded above by the coefficients of the cd-index of a dual stacked polytope of 
dimension $d$ with $d+k+1$ facets. 
The dual stacked polytopes are the polar duals of the stacked polytopes of \cite{Klei-Lee}.
This upper bound would be sharp because the structural recursion can be  used to construct Bruhat intervals which are the face lattices of 
duals of stacked polytopes (Proposition~\ref{dual stacked}).

Stanley \cite{Stanley2} conjectured the non-negativity of the cd-indices of a much more general class of Eulerian posets.
We show (Theorem~\ref{lower}) that if the conjectured non-negativity holds for Bruhat intervals, then the cd-index of any lower Bruhat 
interval is bounded above by the cd-index of a Boolean algebra.
Since the flag h-vectors of Bruhat intervals are non-negative, we are able to prove that the flag h-vectors of lower Bruhat intervals are 
bounded above by the flag h-vectors of Boolean algebras (Theorem~\ref{boolean h}).

The remainder of the paper is organized as follows:
We begin with background information on the basic objects appearing in this paper, namely, posets, Coxeter groups, Bruhat order and polytopes in
Section~\ref{pre}, CW complexes and PL topology in Section~\ref{CW} and the cd-index in Section~\ref{cd}.
In Section \ref{zipping}, the zipping operation is introduced, and its basic properties are proven.
Section \ref{Building} contains the proof of the structural recursion.
In Section \ref{recursion section} we state and prove the cd-index recursions and apply them to determine the 
affine span of cd-indices of Bruhat intervals.
Section \ref{bounds section} is a discussion of conjectured bounds on the coefficients of the cd-index of a Bruhat interval, 
including the construction of Bruhat intervals which are isomorphic to the face lattices of  dual stacked polytopes.

\section{Preliminaries}
\label{pre}
In this section we give background information on posets, Coxeter groups, Bruhat order and polytopes.

\subsection*{Posets}
The poset terminology and notation here generally agree with \cite{Stanley1}.
Throughout this paper, all posets considered are finite.

Let $P$ be a poset.
Given $x,y\in P$, say $x$ {\em covers} $y$ and write ``$x\covers y$'' if $x>y$ and if there is no $z\in P$ with $x>z>y$.
Given $x\in P$, define $D(x):=\set{y\in P : y< x}$. 
If $P$ has a unique minimal element, it is denoted $\0$, and if there is a unique maximal element, it is called $\1$.
The dual of $P$ is the same set of elements with the reversed partial order.
An {\em induced subposet} of $P$ is a subset $S\subseteq P$, together with the partial order induced on $S$ by $\le_P$.
Often this is referred to simply as a {\em subposet}.

A {\em chain} is a totally ordered subposet of $P$.
A chain is {\em maximal} if it is not properly contained in any other chain.
A poset is {\em graded} if every maximal chain has the same number of elements.
A {\em rank} function on a graded poset $P$ is the unique function such that $\rank(x)=0$ for any minimal element $x$, and 
$\rank(x)=\rank(y)+1$ if $x\covers y$.
Given two posets $P$ and $Q$, form their {\em product}  $P\times Q$.
The underlying set is the ordered pairs $(p,q)$ with $p\in P$ and $q\in Q$, and the partial order is $(p,q)\le_{P\times Q} (p',q')$ if and only if 
$p\le_P p'$ and $q\le_Q q'$.
The product of $P$ with a two-element chain is called the pyramid $\Pyr(P)$.
A poset $Q$ is an {\em extension} of $P$ if the two are equal as sets, and if $a\le_Pb$ implies $a\le_Qb$.
The {\em join} $x\join y$ of two elements $x$ and $y$ is the unique minimal element in $\set{z:z\ge x,z\ge y}$, if it exists.
The {\em meet} $x\meet y$ is the  unique maximal element in $\set{z:z\le x,z\le y}$ if it exists.
A poset is called a {\em lattice} every pair of elements $x$ and $y$ has a meet and a join.

Given posets $P$ and $Q$, a map $\eta:P\rightarrow Q$ is {\em order-preserving} if $a\le_Pb$ implies $\eta(a)\le_Q\eta(b)$.
Consider the set $\bar{P}:=\set{\eta^{-1}(q):q\in Q}$ of fibers of $\eta$, and define a relation $\le_{\bar{P}}$ on $\bar{P}$ by 
$F_1 \le_{\bar{P}}F_2$ 
if there exist $a \in F_1$ and $b \in F_2$ such that $a \le_P b$.
If $\le_{\bar{P}}$ is a partial order, $\bar{P}$ is called the {\em fiber poset} of $P$ with respect to $\eta$.
In this case, there is a surjective order-preserving map $\nu:P \rightarrow \bar{P}$ given by 
        $\nu:a \mapsto \eta^{-1}(\eta(a))$, and an injective order-preserving map 
        $\bar{\eta}:\bar{P} \rightarrow Q$ such that $\eta = \bar{\eta}\circ \nu$.
Call $\eta$ an {\em order-projection} if it is order-preserving and has the following property: 
For all $q \le r$ in $Q$, there exist $a \le b \in P$ with $\eta(a)=q$ and $\eta(b)=r$.
In particular, an order-projection is surjective.

\begin{prop}
\label{o-pproperties}
Let $\eta:P\rightarrow Q$ be an order-projection.
Then 
\begin{enumerate}
\item[(i) ] $\bar{P}$ is a fiber poset.
\item[(ii) ] $\bar{\eta}$ is an order-isomorphism.
\end{enumerate}
\end{prop}
\begin{proof}
Assertion (i) is just the statement that $\le_{\bar{P}}$ is a partial order.
The reflexive property is trivial.
Let $A = \eta^{-1}(q)$ and $B = \eta^{-1}(r)$ for $q,r \in Q$.
If $A \le_{\bar{P}} B$ and $B \le_{\bar{P}} A$, we can find $a_1,a_2,b_1,b_2$ with $\eta(a_1)=\eta(a_2)=q$,
        $\eta(b_1)=\eta(b_2)=r$, $a_1 \le b_1$ and $b_2 \le a_2$.
Because $\eta$ is order-preserving, $q \le r$ and $r \le q$, so $q=r$ and therefore $A=B$.
Thus the relation is anti-symmetric.

To show that $\le_{\bar{P}}$ is transitive, suppose $A \le_{\bar{P}} B$ and $B \le_{\bar{P}} C$.
Then there exist $a \in A$, $b_1,b_2\in B$ and $c\in C$ with $\eta(a)=q$, $\eta(b_1)=\eta(b_2)=r$ and $\eta(c)=s$ 
        such that $a \le b_1$ and $b_2 \le c$.
Because $\eta$ is order-preserving, we have $q \le r \le s$.
Because $\eta$ is an order-projection, one can find $a' \le c' \in P$ with $\eta(a')=q$ and $\eta(c')=s$.
So $A \le_{\bar{P}} C$.

Since $\eta$ is surjective, $\bar{\eta}$ is an order-preserving bijection.
Let $q \le r$ in $Q$.
Then, because $\eta$ is an order-projection, there exist $a \le b \in P$ with $\eta(a)=q$ and $\eta(b)=r$.
So $\bar{\eta}^{-1}(q) = \eta^{-1}(q) \le \eta^{-1}(r) = \bar{\eta}^{-1}(r)$ in $\bar{P}$.
Thus $\bar{\eta}^{-1}$ is order-preserving.
\end{proof}

\subsection*{Coxeter groups and Bruhat order}
Here we give the definition of a Coxeter group and the Bruhat order (or strong order) on a Coxeter group, as well as one alternate 
characterization of the Bruhat order.
Further information on Coxeter groups and the Bruhat order can be found for example in \cite{Bourbaki,Humphreys}.

A {\em Coxeter system} is a pair $(W,S)$, where $W$ is a group, $S$ is a set of generators, and $W$ is given by 
        the presentation $(st)^{m(s,t)} = 1$ for all $s,t \in S$, with the requirements that:
\begin{enumerate}
\item[(i) ] $m(s,s) = 1$ for all $s \in S$, and 
\item[(ii) ] $2 \le m(s,t) \le \infty$ for all $s \ne t$ in $S$.
\end{enumerate}
We use the convention that $x^\infty=1$ for any $x$, so that $(st)^\infty=1$ is a trivial relation.
The Coxeter system is called  {\em universal} if $m(s,t)=\infty$ for all $s\neq t$.
We will refer to a ``Coxeter group,'' $W$ with the understanding that a generating set $S$ has been chosen such that $(W,S)$ is a Coxeter system.
In what follows, $W$ or $(W,S)$ will always refer to a fixed Coxeter system, and $w$ will be an element of $W$.
Examples of finite Coxeter groups include the symmetric group, other Weyl groups of root systems, and symmetry groups
        of regular polytopes.

Readers not familiar with Coxeter groups should concentrate on the symmetric group $S_n$ of permutations of the numbers $\set{1,2,\ldots,n}$.
In particular, some of the figures will illustrate the case of $S_4$.
Let $r$ be the transposition $(12)$, let $s:=(23)$  and $t:=(34)$.
Then $(S_4,\set{r,s,t})$ is a Coxeter system with $m(r,s)=m(s,t)=3$ and $m(r,t)=2$.

Call a word $w = s_1s_2\cdots s_k$ with letters in $S$ a {\em reduced} word for $w$ if $k$ is as small as possible.
Call this $k$ the {\em length} of $w$, denoted $l(w)$.
We will use the symbol ``1'' to represent the empty word, which corresponds to the identity element of $W$.
Given any words $a_1$ and $a_2$ and given words $b_1=stst\cdots$ with $l(b_1)=m(s,t)$ and $b_2:=tsts\cdots$ with $l(b_2)=m(s,t)$, the words 
$a_1b_1a_2$ and $a_1b_2a_2$ both stand for the same element.
Such an equivalence is called a {\em braid move}.
A theorem of Tits says that given any two reduced words $a$ and $b$ for the same element, $a$ can be transformed into $b$ by a sequence of 
braid moves.

There are several equivalent definitions of Bruhat order.  
See \cite{Deodhar} for a discussion of the equivalent formulations.
One definition is by the ``Subword Property.''
Fix a reduced word $w = s_1s_2\cdots s_k$.
Then $v \le_B w$ if and only if there is a reduced subword $s_{i_1}s_{i_2}\cdots s_{i_j}$ corresponding to $v$ such that 
$1 \le i_1 < i_2 < \cdots < i_j \le k.$
We will write $v \le w$ for $v \le_B w$ when the context is clear.

Bruhat order is ranked by length.
The element $w$ covers the elements which can be represented by reduced words obtained by deleting a single letter from a reduced word for $w$.
We will need the ``lifting property'' of Bruhat order, which can be proven easily using the Subword Property.

\begin{prop}
\label{lifting}
If $w\in W$ and $s\in S$ have $w>ws$ and $us>u$, then the following are equivalent:
\begin{enumerate}
\item[(i) ] $w>u$
\item[(ii) ] $ws>u$
\item[(iii) ] $w>us$
\end{enumerate}
\end{prop}

\subsection*{Polytopes}
A {\em (convex) polytope} $\Po$ is the convex hull of a finite number of points, or equivalently a bounded set which is the intersection of 
a finite number of closed halfspaces.
An {\em affine subspace} of $\reals^d$ is a subset of $\reals^d$ which can be written as $L+v:=\set{x+v:x\in L}$ for some linear subspace $L$ and some vector $v\in \reals^d$.
The {\em affine span} of a set $S\subseteq\reals^d$ is the intersection of all affine subspaces containing $S$.
The {\em dimension} of a polytope $\Po$ is the dimension of its affine span.

A {\em hyperplane} is the set of points $x$ satisfying $a\cdot x=b$ for some fixed $a$ and $b$.
It is called a {\em supporting} hyperplane of a polytope $\Po$ if $a\cdot x\le b$ for every point in $\Po$ or if $a\cdot x\ge b$ for every point in $\Po$. 
A {\em face} of $\Po$ is any intersection of $\Po$ with a supporting hyperplane.
In particular $\emptyset$ is a face of $\Po$, and any face of $\Po$ is itself a convex polytope.
By convention, $\Po$ is also a face of $\Po$.
A {\em facet} of $\Po$ is a face whose dimension is one less than the dimension of $\Po$.
The {\em face lattice} of $\Po$ is the set of faces of $\Po$ partially ordered by inclusion, and this partial order can be shown to be a lattice.
Two polytopes are of the same {\em combinatorial type} if their face lattices are isomorphic as posets.

We will need two geometric constructions on polytopes, the pyramid operation $\Pyr$ and the vertex-shaving operation $\Sh_v$.
Given a polytope $P$ of dimension $d$, $\Pyr(P)$ is the convex hull of the union of $P$ with some vector $v$ which is not in the affine span 
of $\Po$.
This is unique up to combinatorial type and the face poset of $\Pyr(\Po)$ is just the pyramid of the face poset of $\Po$.

Consider a polytope $P$ and a chosen vertex $v$.
Let $H = \set{a\cdot x = b}$ be a hyperplane that separates $v$ from the other vertices of $P$.
In other words, $a\cdot v > b$ and $a\cdot v' <b$ for all vertices $v' \ne v$.
Then the polytope $\Sh_v(P) = P \cap \set{a\cdot x \le b}$ is called the {\em shaving} of $P$ at $v$.
This is unique up to combinatorial type.
Every face of $P$, except $v$, corresponds to a face in $\Sh_v(P)$ and, in addition, for every face of $P$ strictly 
        containing $v$, there is an additional face of one lower dimension in $\Sh_v(P)$.
In Section \ref{CW} we describe how this operator can be extended to regular CW spheres, and in Section \ref{Building} we 
describe the corresponding operator on posets.

Further information on polytopes can be found for example in \cite{Polytopes}.

\section{CW complexes and PL topology}
\label{CW}
This section provides background material on finite CW complexes and PL topology which will be useful in Section \ref{zipping}.
More details about CW complexes, particularly as they relate to posets, can be found in \cite{Bj}. 
Additional details about PL topology can be found in \cite{OrientedMatroids, Rourke-Sanderson}.

A set of $n$ points in $\reals^d$ is {\em affinely independent} if the smallest affine subspace containing them has dimension $n-1$.
A {\em simplex} is a polytope which is the convex hull of an affinely independent set $I$ of points.
The faces of the simplex are the convex hulls of subsets of $I$.
A {\em geometric simplicial complex}~$\Delta$ is a finite collection of simplices (called {\em faces} of the complex) 
such that
\begin{enumerate}
\item[(i) ] If $\sigma\in\Delta$ and $\tau$ is a face of $\sigma$, then $\tau\in\Delta$.
\item[(ii) ] If $\sigma,\tau\in \Delta$ then $\sigma\cap\tau$ is a face of $\sigma$ and of $\tau$.
\end{enumerate}
The zero-dimensional faces are called {\em vertices}.
The {\em underlying space} $|\Delta|$ of~$\Delta$ is the union in $\reals^d$ of the faces of~$\Delta$.

An {\em abstract simplicial complex}~$\Delta$ on a finite {\em vertex-set} $V$ is a collection of subsets of $V$ called {\em faces}, 
with the following properties:
\begin{enumerate}
\item[(i) ]Every singleton is a face.
\item[(ii) ]Any subset of a face is another face.
\end{enumerate}

Given a geometric simplicial complex~$\Delta$ with vertices $V$, the collection
\[\set{F\subseteq V: \mbox{the convex hull of }F\mbox{ is in }\Delta}\] is an abstract simplicial complex.
The process can be reversed: given an abstract simplicial complex~$\Delta$, there is a construction which produces a geometric simplicial complex 
whose underlying abstract simplicial complex is exactly~$\Delta$.
This {\em geometric realization} is unique up to homeomorphism, so it makes sense to talk about the topology of an abstract simplicial complex.
Two geometric simplicial complexes are {\em combinatorially isomorphic} if their underlying abstract simplicial 
complexes are isomorphic.
If two complexes are combinatorially isomorphic then their underlying spaces are homeomorphic, but the converse is not true.

Given simplicial complexes~$\Delta$ and~$\Gamma$, say~$\Gamma$ is {\em subdivision} of~$\Delta$ if $|\Gamma|=|\Delta|$ and if 
every face of $\Gamma$ is contained in some face of~$\Delta$.
A simplicial complex is a {\em PL $d$-sphere} if it admits a simplicial subdivision which is combinatorially isomorphic to some simplicial 
subdivision of the boundary of a $(d+1)$-dimensional simplex.
A simplicial complex is a {\em PL $d$-ball} if it admits a simplicial subdivision which is combinatorially isomorphic to some simplicial 
subdivision of a $d$-dimensional simplex.

We now quote some results about PL balls and spheres.
Some of these results appear topologically obvious but, surprisingly, not all of these statement are true with the ``PL '' deleted.
This is the reason that we introduce PL balls and spheres, rather than dealing with ordinary topological balls and spheres.
Statement (iii) is known as Newman's Theorem.

\begin{theorem}\cite[Theorem 4.7.21]{OrientedMatroids}
\label{PL tricks}
\begin{enumerate}
\item[(i) ]Given two PL $d$-balls whose intersection is a PL $(d-1)$-ball lying in the boundary of each, the union of the two is a PL $d$-ball. 
\item[(ii) ]Given two PL $d$-balls whose intersection is the entire boundary of each, the union of the two is a PL $d$-sphere. 
\item[(iii) ]The closure of the complement of a PL $d$-ball embedded in a PL $d$-sphere is a PL $d$-ball.
\qed
\end{enumerate}
\end{theorem}

Given two abstract simplicial complexes~$\Delta$ and~$\Gamma$, let $\Delta * \Gamma$ be the {\em join} of~$\Delta$ and~$\Gamma$, 
a simplicial complex whose vertex set  is the disjoint union of the vertices of~$\Delta$ and of~$\Gamma$, and whose faces are exactly the 
sets $F\cup G$ for all faces $F$ of $\Delta$ and $G$ of~$\Gamma$.
Let $B^d$ stand for a PL $d$-ball, and let $S^d$ be a PL $d$-sphere.

\begin{prop}\cite[Proposition 2.23]{Rourke-Sanderson}
\label{joins}
\begin{eqnarray*}
  B^p*B^q&\cong&B^{p+q+1}\\
  S^p*B^q&\cong&B^{p+q+1}\\
  S^p*S^q&\cong&S^{p+q+1}
\end{eqnarray*}
\qed
\end{prop}
Here $\cong$ stands for PL homeomorphism, a stronger condition than homeomorphism which requires a compatibility of PL-structures as well.
The point is that $B^p*B^q$ is a PL ball, etc.

Given a poset $P$ the order complex $\Delta(P)$ is the abstract simplicial complex whose vertices are the elements of $P$ and whose faces 
are the chains of $P$.
The order complex of an interval $[x,y]$ will be written $\Delta[x,y]$, rather than $\Delta([x,y])$, and similarly $\Delta(x,y)$ instead of 
$\Delta((x,y))$. 
When $P$ is a poset with a $\0$ and a $\1$, by convention, statements about the topology of $P$ are understood to apply to the order complex of 
$(\0,\1)=P-\set{\0,\1}$.
Thus for example, the statement that ``$P$ is a PL sphere'' means that $\Delta(\0,\1)$ is a PL sphere.
The following proposition follows immediately from \cite[Theorem 4.7.21(iv)]{OrientedMatroids}:

\begin{prop}
\label{interval}
If $P$ is a PL sphere then any interval $[x,y]_P$ is a PL sphere.
\qed
\end{prop}

An {\em open cell} is any topological space isomorphic to an open ball.
A {\em CW complex} $\Omega$ is a Hausdorff topological space with a decomposition as a disjoint union of cells, such that
for each cell $e$, the homeomorphism mapping an open ball to $e$ is required to extend to a continuous map from the closed ball to $\Omega$.
The image of this extended map is called a {\em closed cell}, specifically the {\em closure} of $e$. 
The {\em face poset} of $\Omega$ is the set of closed cells, together with the empty set, partially ordered by containment.
The {\em $k$-skeleton} of $\Omega$ is the union of the closed cells of dimension $k$ or less.
A CW complex is  {\em regular} if all the closed cells are homeomorphic to closed balls.

Call $P$ a CW poset if it is the face poset of a regular CW complex $\Omega$.
It is well known that in this case $\Omega$ is homeomorphic to $\Delta(P-\set{\0})$.
The following theorem is due to Bj\"{o}rner~\cite{Bj}.
\begin{theorem}
\label{CW thm}
A non-trivial poset $P$ is a CW poset if and only if
\begin{enumerate}
\item[(i) ] $P$ has a minimal element $\hat{0}$, and
\item[(ii) ] For all $x\in P-\set{\hat{0}}$, the interval $[\hat{0},x]$ is a sphere.
\qed
\end{enumerate}
\end{theorem}
Given a CW poset, Bj\"{o}rner constructs a complex $\Omega(P)$ recursively by constructing the $(k-1)$-skeleton, and then attaching 
$k$-cells in a way that agrees with the order relations in $P$.

The polytope operations $\Pyr$ and $\Sh_v$ can also be defined on regular CW spheres.
Both operations preserve PL sphericity by Theorem~\ref{PL tricks}(ii).
We give informal descriptions which are easily made rigorous.
Consider a regular CW $d$-sphere $\Omega$ embedded as the unit sphere in $\reals^{d+1}$.
The new vertex in the $\Pyr$ operation will be the origin.
Each face of $\Omega$ is also a face of $\Pyr(\Omega)$ and for each nonempty face $F$ of $\Omega$ there is a new face $F'$ of 
$\Pyr(\Omega)$, described by 
\[F':=\set{v\in \reals^{d+1}:0<|v|<1,\frac{v}{|v|}\in F}.\]
The set $\set{v\in\reals^{d+1}:|v|>1}\cup\set{\infty}$ is also a face of $\Pyr(\Omega)$ (the ``base'' of the pyramid) where $\infty$ is the point 
at infinity which makes $\reals^{d+1}\cup\set{\infty}$ a $(d+2)$-sphere.

Consider a regular CW sphere $\Omega$ and a chosen vertex $v$.
Adjoin a new open cell to make $\Omega'$, a ball of one higher dimension.
Choose $S$ to be a small sphere $|x-v|=\epsilon$, such that the only vertex inside the sphere is $v$ and the only faces which intersect $S$ are faces 
which contain $v$.
(Assuming some nice embedding of $\Omega$ in space, this can be done.)
Then $\Sh_v(\Omega)$ is the boundary of the ball obtained by intersecting $\Omega'$ with the set $|x-v|\ge \epsilon$.
As in the polytope case, this is unique up to combinatorial type.
Every face of $\Omega$, except $v$, corresponds to a face in $\Sh_v(\Omega)$, 
and for every face of $\Omega$ strictly containing $v$, there is an additional face of one lower dimension in $\Sh_v(\Omega)$.

Given a poset $P$ with $\0$ and $\1$, call $P$ a regular CW sphere if $P-\set{\hat{1}}$ is the face poset of a regular CW complex which is a sphere.
By Theorem~\ref{CW thm}, $P$ is a regular CW sphere if and only if every lower interval of $P$ is a sphere.
In light of Proposition~\ref{interval}, if $P$ is a PL sphere, then it is also a CW sphere, but not conversely.
Section \ref{Building} describes a construction on posets which corresponds to $\Sh_v$.

\section[cd-Index]{The cd-index of an Eulerian poset}
\label{cd}
In this section we give the definition of Eulerian posets, flag f-vectors, flag h-vectors, and the cd-index, and quote results about the cd-indices 
of polytopes.

The {\em M\"obius function} $\mu\colon \set{(x,y)\colon x\le y \mbox{ in } P} \rightarrow \integers$ is defined recursively by setting $\mu(x,x) = 1$
for all $x\in P$, and 
        \[\mu(x,y) = -\sum_{x\le z <y} \mu(x,z) \mbox{ for all } x < y \mbox{ in } P.\]
A poset $P$ is {\em Eulerian} if $\mu(x,y) = (-1)^{\rank(y)-\rank(x)}$ for all intervals $[x,y] \subseteq P$.
For a survey of Eulerian posets, see \cite{Stanley3}.

Verma \cite{Verma} gives an inductive proof that Bruhat order is Eulerian, by counting elements of even and odd rank.
Rota \cite{Rota} proved that the face lattice of a convex polytope is an Eulerian poset  (See also \cite{Lindstr}). 
More generally, the face poset of a CW sphere is Eulerian.
In \cite{Bj}, Bj\"{o}rner showed that Bruhat intervals are CW spheres.

Let $P$ be a  graded poset, rank $n+1$, with a minimal element $\hat{0}$ and a maximal element $\hat{1}$.
For a chain $C$ in $P-\set{\0,\1}$, define $\rank(C) = \set{\rank(x):x\in C}$.
Let $[n]$ denote the set of integers $\set{1,2,\ldots,n}$.
For any $S \subseteq [n]$, define
        \[\alpha_P(S)= \# \set{\mbox{chains } C \subseteq P:\rank(C) = S}.  \]
The function $\alpha_P: 2^{[n]}\rightarrow \naturals$ is called the {\em flag f-vector}, because it is a refinement of 
        the {\em f-vector}, which counts the number of elements of each rank.

Define a function $\beta_P: 2^{[n]}\rightarrow \naturals$ by
\[                \beta_P(S) = \sum_{T\subseteq S} (-1)^{|S-T|}\alpha_P(T).\]
The function $\beta_P$ is called the {\em flag h-vector} of $P$ because of its relation to the usual {\em h-vector}.

Bayer and Billera \cite{Bayer-Billera} proved a set of linear relations on the flag f-vector of an Eulerian poset, called the 
Generalized Dehn-Sommerville relations.
They also proved that the Generalized Dehn-Sommerville relations and the relation $\alpha_P(\emptyset) = 1$ are the complete set of 
affine relations satisfied by flag f-vectors of all Eulerian posets.

Let $\integers \langle a,b \rangle$ be the vector space of {\em ab-polynomials}, that is, polynomials over non-commuting variables 
$a$ and $b$ with integer coefficients.
Subsets $S \subseteq [n]$  can be represented by monomials $u_S = u_1u_2\cdots u_n \in \integers \langle a,b \rangle$,
        where $u_i = b$ if $i \in S$ and $u_i = a$ otherwise.
Define ab-polynomials $\Upsilon_P$ and $\Psi_P$ to encode the flag f-vector and flag h-vector respectively.
\begin{eqnarray*}
        \Upsilon_P(a,b) &:=& \sum_{S\subseteq [n]} \alpha_P(S)u_S\\
        \Psi_P(a,b) &:=& \sum_{S\subseteq [n]} \beta_P(S)u_S.
\end{eqnarray*}
The polynomial $\Psi_P$ is commonly called the {\em ab-index}.
There is no standard name for $\Upsilon_P$, but here we will call it the {\em flag index}.
It is easy to show that $\Upsilon_P(a-b,b) = \Psi_P(a,b)$.

Let $c = a+b$ and $d = ab + ba$ in $\integers \langle a,b \rangle$.
The flag f-vector of a graded poset $P$ satisfies the Generalized Dehn-Sommerville relations if and only if $\Psi_P(a,b)$ can be 
        written as a polynomial in $c$ and $d$ with integer coefficients, called the {\em cd-index} of $P$.
This surprising fact was conjectured by J.\ Fine and proven by Bayer and Klapper \cite{Bayer-Klapper}.
The cd-index is {\em monic}, meaning that the coefficient of $c^n$ is always~1.
The existence and monicity of the cd-index constitute the complete set of affine relations on the flag f-vector of an Eulerian poset.
Setting the degree of $c$ to be 1 and the degree of $d$ to be 2, the cd-index of a poset of rank $n+1$ is homogeneous of degree $n$.
The number of cd-monomials of degree $n-1$ is $F_n$, the $n^{\mbox{\small th}}$ Fibonacci number, with $F_1 = F_2 = 1$.
Thus the affine span of flag f-vectors of Eulerian posets of degree $n$ has dimension $F_n -1$.

The literature is divided on notation for the cd-index, due to two valid points of view as to what the ab-index is.
If one considers $\Psi_P$ to be a polynomial function of non-commuting variables $a$ and $b$, one may consider the cd-index to be a 
        different polynomial function in $c$ and $d$, and give it a different name, typically $\Phi_P$.
On the other hand, if $\Psi_P$ is a vector in a space of ab-polynomials, the cd-index is the same vector,
        which happens to be written as a linear combination of monomials in $c$ and~$d$.
Thus one would call the cd-index $\Psi_P$.
We will primarily use the notation $\Psi_P$, except that when we talk about inequalities on the coefficients of the cd-index, we use $\Phi_P$.

Aside from the existence and monicity of the cd-index, there are no additional affine relations on flag f-vectors of polytopes.  
Bayer and Billera \cite{Bayer-Billera} and later Kalai \cite{Kalai} gave a basis of polytopes whose flag f-vectors 
        span $\integers \langle c,d \rangle$.
Much is also known about bounds on the coefficients of the cd-index of a polytope.
A bound on the cd-index implies bounds on $\alpha$ and $\beta$, because $\alpha$ and $\beta$ can be written as positive 
        combinations of coefficients of the cd-index.
The first consideration is the non-negativity of the coefficients.
Stanley \cite{Stanley2} conjectured that the coefficients of the cd-index are non-negative whenever $P$ triangulates a homology sphere 
(or in other words when $P$ is a {\em Gorenstein* poset}).
He also showed that the coefficients of $\Phi_P$ are non-negative for a class of CW-spheres which includes convex polytopes.

Ehrenborg and Readdy described how the cd-index is changed by the poset operations of pyramid and vertex shaving.
The following is a combination of Propositions 4.2 and 6.1 of \cite{Coproducts}.

\begin{prop}
\label{cd ops}
Let $P$ be a graded poset and let $a$ be an atom.  Then
\begin{eqnarray*}
\Psi_{\Pyr(P)} & = & \frac{1}{2}\left(\Psi_P\cdot c + c\cdot\Psi_P + \sum_{x\in P,\hat{0}<x<\hat{1}}\Psi_{[\hat{0},x]}\cdot d \cdot \Psi_{[x,\hat{1}]} \right)\\
\Psi_{\Sh_a(P)} &=& \Psi_P + \frac{1}{2}\left(\Psi_P\cdot c - c\cdot \Psi_P + \sum_{a<x<\1}\Psi_{[a,x]}\cdot d\cdot\Psi_{[x,\1]} \right).
\end{eqnarray*}
\end{prop}

Ehrenborg and Readdy also defined a derivation on cd-indices and used it to restate the formulas in Proposition \ref{cd ops}.
The derivation $G$ (called $G'$ in \cite{Coproducts}) is defined by $G(c)=d$ and $G(d)=dc$.
The following is a combination of Theorem 5.2 and Proposition 6.1 of \cite{Coproducts}.

\begin{prop}
\label{cd derivation}
Let $P$ be a graded poset and let $a$ be an atom.  Then
\begin{eqnarray*}
\Psi_{\Pyr(P)} & = & c \cdot \Psi_P + G\left(\Psi_P\right)\\
\Psi_{\Sh_a(P)} &=& \Psi_P + G\left(\Psi_{[a,\1]}\right).
\end{eqnarray*}
\end{prop}

\begin{cor}
\label{Pyr}
Let $P$ be a homogeneous cd-polynomial whose lexicographically first term is $T$.
Then the lexicographically first term of $\Pyr(P)$ is $c\cdot T$.
In particular, the kernel of the pyramid operation is the zero polynomial.
\qed
\end{cor}

\section{Zipping}
\label{zipping}
In this section we introduce the zipping operation and prove some of its important properties.
In particular, zipping will be part of a new inductive proof that Bruhat intervals are spheres and thus Eulerian.
A {\em zipper} in a poset $P$ is a triple of distinct elements $x,y,z \in P$ with the following properties:
\begin{enumerate}
\item[(i) ]$z$ covers $x$ and $y$ but covers no other element.
\item[(ii) ]$z = x \vee y$.
\item[(iii) ]$D(x)=D(y)$.
\end{enumerate} 
Call the zipper {\em proper} if $z$ is not a maximal element.
If $(x,y,z)$ is a zipper in $P$ and $[a,b]$ is an interval in $P$ with $x,y,z\in [a,b]$ then $(x,y,z)$ is a zipper in $[a,b]$.

Given $P$ and a zipper $(x,y,z)$ one can ``zip'' the zipper as follows:
Let $xy$ stand for a single new element not in $P$.
Define $P' = (P-\set{x,y,z}) \cup \set{xy}$, with a binary relation called $\preceq$, given by:
        \begin{eqnarray*}
                & a  \preceq b    & \mbox{ if } a\le b\\
                & xy \preceq a    & \mbox{ if } x\le a \mbox{ or if } y \le a \\
                & a  \preceq xy   & \mbox{ if } a\le x \mbox{ or (equivalently) if } a \le y \\
                & xy \preceq xy   & \end{eqnarray*}

For convenience, $[a,b]$ will always mean the interval $[a,b]_\le$ in $P$ and $[a,b]_\preceq$ will mean an
        interval in $P'$.
In each of the following propositions, $P'$ is obtained from $P$ by zipping the proper zipper $(x,y,z)$, although some of the results are true 
even when the zipper in not proper.

\begin{prop}
\label{poset}
$P'$ is a poset under the partial order $\preceq$.
\end{prop} 

\begin{proof}
One sees immediately that $\preceq$ is reflexive and that antisymmetry holds in $P'-\set{xy}$.
If $xy \preceq a$ and $a \preceq xy$, but $a\ne xy$, then $a\in P-\set{x,y,z}$.
We have $a\le x$ and $a\le y$.
Also, either $x\le a$ or $y\le a$.
By antisymmetry in $P$, either $a=x$ or $a=y$.
This contradiction shows that $a=xy$.
Transitivity follows immediately from the transitivity of $P$ except perhaps when $a\preceq xy$ and $xy\preceq b$.
In this case, $a\le x$ and $a\le y$.
Also, either $x\le b$ or $y\le b$.
In either case, $a\le b$ and therefore $a\preceq b$.
\end{proof}

\begin{prop}
\label{mu1}
If $a\preceq xy$ then $\mu_{P'}(a,xy)=\mu_P(a,x)=\mu_P(a,y)$.
If $a\preceq b \in P'$ with $a\ne xy$, then $\mu_{P'}(a,b)=\mu_P(a,b)$.
\end{prop} 

Suppose $[a,b]_\preceq$ is any non-trivial interval in $P'$.
If $a\not\preceq xy$, then $[a,b]_\preceq=[a,b]_\le$.
If $b=xy$, then $[a,b]_\preceq\cong[a,x]_\le$.
If $b\ne xy$ and $b\not>z$, then $[a,b]_\le$ does not contain both $x$ and $y$, and we obtain $[a,b]_\preceq$ from $[a,b]_\le$ by replacing $x$ or $y$ by 
$xy$ if necessary.
Thus in the proofs that follow, one needs only to check two cases: the case where $a\prec xy$ and $b>z$ and the case where $a=xy$.

\begin{proof}[Proof of Proposition \ref{mu1}]
Let $a\preceq b$ with $a\ne xy$.
One needs only to check the case where $a\prec xy$ and $b>z$.
This is done by induction on the length of the longest chain from $z$ to $b$.
If $b\covers z$ then 
\begin{eqnarray*}
\mu_P(a,b) & = & -\mu_P(a,z)-\mu_P(a,x)-\mu_P(a,y)-\sum_{a\le p < b: p\ne x,y,z}\mu_P(a,p) \\
           & = & \sum_{a\le p < x}\mu_P(a,p) -  \sum_{a\le p < b: p\ne x,y,z}\mu_P(a,p) \\
           & = & -\mu_{P'}(a,xy) - \sum_{a\preceq p \prec b: p\ne xy}\mu_{P'}(a,p) \\
           & = & \mu_{P'}(a,b).
\end{eqnarray*}
Here the second line is obtained by properties (i) and (iii).
If $b$ does not cover $z$, use the same calculation, employing induction to go from the second line to the third line.
\end{proof}

\begin{prop}
\label{mu2}
If $xy\preceq b \in P'$, then $\mu_{P'}(xy,b)=\mu_P(x,b) + \mu_P(y,b) + \mu_P(z,b)$.
\end{prop} 

\begin{proof}
In light of Proposition \ref{mu1}, one can write:
\begin{eqnarray*}
\mu_{P'}(xy,b)  & = & -\sum_{xy\prec p \preceq b}\mu_{P'}(p,b) \\
                & = &  \mu_P(z,b)-\sum_{\substack{p: x< p \le b\\\mbox{{\small or }}\\ y< p \le b}}\mu_P(p,b) \\
                & = &  \mu_P(z,b) -\sum_{x< p \le b}\mu_P(p,b)-\sum_{y< p \le b}\mu_P(p,b) + \sum_{z \le p \le b}\mu_P(p,b)\\
                & = & \mu_P(x,b) + \mu_P(y,b) + \mu_P(z,b).
\end{eqnarray*}
\end{proof}

The following two corollaries follow trivially from Propositions \ref{mu1} and \ref{mu2} and the observation that if $P$ is ranked, then $P'$ inherits a 
rank function.

\begin{cor}
\label{thin}
If $P$ is thin, then so is $P'$
\qed
\end{cor} 

\begin{cor}
\label{eul}
If $P$ is Eulerian, then so is $P'$.
\qed
\end{cor} 

\begin{theorem}
\label{cd thm}
$P$ has a cd-index if and only if $P'$ has cd-index.  
The cd-indices are related by: 
        \[\Psi_{P'}=\Psi_P-\Psi_{[\0,x]_\le}\cdot d \cdot\Psi_{[z,\1]_\le}.\]
\end{theorem} 

\begin{proof}
We subtract from $\Upsilon_P$ the chains which disappear under the zipping.
First subtract the terms which came from chains through $x$ and $z$.
Any such chain is a chain in $[\hat{0},x]_P$ concatenated with a chain in $[z,\hat{1}]_P$.
So the terms subtracted off are $\Upsilon_{[\hat{0},x]_P} \cdot  b \cdot b \cdot \Upsilon_{[z,\hat{1}]}$.
Then subtract a similar term for chains through $y$ and $z$.
In fact, by condition (iii) of the definition of a zipper, the term for chains through $y$ and $z$ is identical to the term for chains through $x$ and $z$.
Subtract $\Upsilon_{[\hat{0},x]_P} \cdot  a \cdot b \cdot \Upsilon_{[z,\hat{1}]}$ for the chains which go through $z$ but skip the rank below $z$.
Finally, $x$ is identified with $y$, so there is a double-count which must be subtracted off.
If two chains are identical except that one goes through $x$ and the other goes through $y$, then they are counted twice in $P$ but only once in $P'$.
Because $x\join y= z$, if such a pair of chains include an element whose rank is $\rank(z)$, then that element is $z$.
But the chains through $z$ have already been subtracted, so we need to subtract off $\Upsilon_{[\hat{0},x]_P} b\cdot a \cdot \Upsilon_{[z,\hat{1}]}$.
We have again used condition (iii) here.
Thus:
\begin{eqnarray*}
     \Upsilon_{P'}&=&\Upsilon_P-\Upsilon_{[\hat{0},x]_P}(2bb+ab+ba)\cdot\Upsilon_{[z,\hat{1}]_P}.
\end{eqnarray*}
Replacing $a$ by $a-b$ one obtains:
\begin{eqnarray}
\label{ab formula}
        \Psi_{P'} & = & \Psi_P - \Psi_{[\hat{0},x]_P}\cdot (ab + ba)\cdot \Psi_{[z,\hat{1}]_P}\\
                  & = & \Psi_P - \Psi_{[\hat{0},x]_P}\cdot d\cdot\Psi_{[z,\hat{1}]_P}. \nonumber
\end{eqnarray}
\end{proof}

\begin{theorem}
\label{pl}
If $P$ is a PL sphere, then so is $P'$.
\end{theorem} 

\begin{proof}
To avoid tedious repetition, we will omit ``PL'' throughout the proof.  All spheres and balls are assumed to be PL.

Suppose $P$ is a $k$-sphere.
Let $\Delta_{xyz}\subset \Delta(\0,\1)$ be the simplicial complex whose facets are maximal chains in $P-\set{\0,\1}$ passing through $x$, $y$ or $z$.
Our first goal is to prove that $\Delta_{xyz}$ is a ball.
Let $\Delta_x\subset \Delta(\0,\1)$ be the simplicial complex whose facets are maximal chains in $(\0,\1)$
        through $x$.
Similarly $\Delta_y$.
One can think of $\Delta_x$ as $\Delta(\0,x)*x*\Delta(x,\1)$.
Thus, by Proposition \ref{joins}, $\Delta_x$ is a $k$-ball, and similarly, $\Delta_y$.
Let $\Gamma = \Delta_x \cap \Delta_y$.
Then~$\Gamma$ is the complex whose facets are almost-maximal chains that can be completed to maximal
        chains either by adding $x$ or $y$.
These are the chains through $z$ which have elements at every rank except at the rank of $x$.
Thus~$\Gamma$ is $\Delta(\0,x)*z*\Delta(z,\1)$, a $(k-1)$-ball, and~$\Gamma$ lies in the boundary of $\Delta_x$, because there is exactly one way to 
complete a facet of~$\Gamma$ to a facet of $\Delta_x$, namely by adjoining~$x$.
Similarly,~$\Gamma$ lies in the boundary of $\Delta_y$.
So by Theorem~\ref{PL tricks}(i), $\Delta_{xyz}=\Delta_x\cup \Delta_y$ is a $k$-ball.

Consider $\Delta((\0,\1)-\set{x,y,z})$, which is the closure of $\Delta(\0,\1)-\Delta_{xyz}$.
By Theorem~\ref{PL tricks}(iii), $\Delta((\0,\1)-\set{x,y,z})$ is also a $k$-ball.
Also consider $\Delta((\0,\1)_\preceq-\set{xy})$, which is isomorphic to $\Delta((\0,\1)-\set{x,y,z})$.
The boundary of $\Delta((\0,\1)_\preceq-\set{xy})$ is a complex whose facets are chains $c$ with the property 
        that for each $c$ there is a unique element of $(\0,\1)_\preceq-\set{xy}$ that completes $c$ to a 
        maximal chain.
However, since $(\0,\1)_\preceq$ is thin by Corollary~\ref{thin}, it has the property that any chain of length $k-1$ can be 
        completed to a maximal chain in $(\0,\1)_\preceq$ in exactly two ways.
Therefore every facet of the boundary of $\Delta((\0,\1)_\preceq-\set{xy})$ is contained in a chain through $xy$.
So $\Delta((\0,\1)_\preceq)$ is the union of a $k$-ball $\Delta((\0,\1)_\preceq-\set{xy})$ with the pyramid over the 
        boundary of $\Delta((\0,\1)_\preceq-\set{xy})$.
By Theorem~\ref{PL tricks}(ii), $\Delta((\0,\1)_\preceq)$ is a $k$-sphere.
\end{proof}

In the case where $P$ is thin, the conditions for a zipper can be simplified.

\begin{prop}
\label{simplezipper}
If $P$ is thin, then (i) implies (iii).
Thus $(x,y,z)$ is a zipper if and only if it satisfies conditions (i) and (ii).
\end{prop}
\begin{proof}
Suppose condition (i) holds but $[\hat{0},x)\ne[\hat{0},y)$.
Then without loss of generality $x$ covers some $a$ which $y$ does not cover.
Since $z$ covers no element besides $x$ and $y$, $[a,z]$ is a chain of length 2, contradicting thinness.
\end{proof}

\section[Building intervals]{Building intervals in Bruhat order}
\label{Building}
In this section we state and prove the structural recursion for Bruhat intervals.
When $s \in S$, $u<us$ and $w<ws$, define a map $\eta:[u,w]\times [1,s]\rightarrow [u,ws]$, as follows:
        \begin{eqnarray*}
                \eta(v,1) & = & v\\
                \eta(v,s) & = & \left\{\begin{array}{ll}
                                        vs & \mbox{if } vs > v\\
                                        v  & \mbox{if } vs < v. \end{array} \right. \end{eqnarray*}
To show that $\eta$ is well-defined, let $v \in [u,w]$.
Then $\eta(v,1) = v \in [u,ws]$ because $ws>w\ge v \ge u$.
Either $\eta(v,s) = v \in [u,ws]$ or $\eta(v,s) = vs$.
In the latter case, $vs>v$, so $u < vs \le ws$ by the lifting property.

\begin{proposition}
\label{eta}
If $u <us$ and $w < ws$, then $\eta:[u,w] \times [1,s] \rightarrow [u,ws]$ is an order-projection.
\end{proposition}
\begin{proof}
To check that $\eta$ is order-preserving, suppose $(v_1,a_1) \le (v_2,a_2)$ in $[u,w] \times [1,s]$.
We have to break up into cases to check that $\eta(v_1,a_1) \le \eta(v_2,a_2)$.
\begin{enumerate}
\item[]
\begin{enumerate}
\item[Case 1:] $a_1 = 1$.\\
If $a_2 = 1$ as well, then $\eta(v_1,a_1) = v_1 \le v_2 =\eta(v_2,a_2) $.
If $a_2 = s$, then $\eta(v_2,a_2)$ is either $v_2$ with $v_2\ge v$ or it is $v_2s$ with $v_2s > v_2 \ge v_1$. 
\item[Case 2:] $a_1 = s$.\\
So $\eta(v_1,a_1)$ is either $v_1$, with $v_1>v_1s$ or it is $v_1s$ with $ v_1s> v_1$. 
We must also have $a_2 = s$, so $\eta(v_2,a_2)$ is either $v_2$ with $v_2>v_2s$ or it is $v_2s$ with $v_2s> v_2$.
If $\eta(v_1,a_1) = v_1$ then $\eta(v_1,a_1) \le v_2 \le \eta(v_2,a_2)$.
If $\eta(v_1,a_1) = v_1s$ and $\eta(v_2,a_2) = v_2$, we have $v_1s>v_1$ and $v_2>v_2s$.
By hypothesis, $v_1 < v_2$, so by the lifting property $\eta(v_1,a_1) = v_1s \le v_2 = \eta(v_1,a_1)$.
If $\eta(v_1,a_1) = v_1s$ and $\eta(v_2,a_2) = v_2s$ then $v_1s>v_1$ and $v_2s>v_2$, so by the lifting property $v_1s\le v_2s$.
\end{enumerate}
\end{enumerate}

It will be useful to identify the inverse image of an element $v \in [u,ws]$.
The inverse image is:
        \[\eta^{-1}(v) = \left\{\begin{array}{ll}
                        \set{(v,1)}              & \mbox{if } v < vs\\
                        \set{(v,1),(vs,s),(v,s)} & \mbox{if } v > vs, \end{array} \right. \]
provided that these elements are actually in $[u,w] \times [1,s]$.
In the case where $v<vs$, we have $u \le v \le w$, where the second inequality is by the 
        lifting property.
So $(v,1)$ is indeed an element of $[u,w] \times [1,s]$.
In the case where $vs < v$, we have by hypothesis $us > u$, so by lifting, $vs \ge u$.
Also by lifting, since $ws>w$ and $v>vs$, we have $w>vs$.
So $(vs,s) \in [u,w] \times [1,s]$.

Now, suppose $x_1 \le x_2 \in [u,ws]$.
To finish the proof that $\eta$ is an order-projection, we must find elements $(v_1,a_1)\le (v_2,a_2)\in [u,w]\times[1,s]$ with $\eta(v_1,a_1)=x_1$ and
        $\eta(v_2,a_2)=x_2$.
Consider 4 cases:
\begin{enumerate}
\item[]
\begin{enumerate}
\item[Case 1:] $x_1<x_1s$ and $x_2<x_2s$.\\
By the inverse-image argument of the previous paragraph, $x_1,x_2 \in [u,w]$, so $\eta(x_1,1)=x_1$, $\eta(x_2,1)=x_2$ and $(x_1,1)\le (x_2,1)$.
\item[Case 2:] $x_1<x_1s$ and $x_2>x_2s$.\\
By the previous paragraph, $x_1,x_2s \in [u,w]$.
Again we have $\eta(x_1,1)=x_1$, and  $\eta(x_2s,s)=x_2$.
By lifting, $x_1\le x_2s$, so $(x_1,1)\le (x_2s,s)$.
\item[Case 3:] $x_1>x_1s$ and $x_2>x_2s$.\\
We have $x_1s,x_2s \in [u,w]$, $\eta(x_1s,s)=x_1$ and $\eta(x_2s,s)=x_2$.
By lifting, $x_1s \le x_2s$, so $(x_1s,s)\le (x_2s,s)$.
\item[Case 4:] $x_1>x_1s$ and $x_2<x_2s$.\\
We have $x_2 \in [u,w]$.
Since $u \le x_1 \le x_2$, $x_1 \in [u,w]$ as well.
So $\eta(x_1,1)=x_1$, $\eta(x_2,1)=x_2$ and $(x_1,1)\le (x_2,1)$.
\end{enumerate}
\end{enumerate}
\end{proof}

In light of the previous section, $\eta$ induces an isomorphism $\bar{\eta}$ between $[u,ws]$ and
        a poset derived from $[u,w] \times [1,s]$, as follows:
For every $v \in [u,w]$ with $vs < v$, ``identify'' $(v,1)$, $(vs,s)$ and $(v,s)$ to make a single element.
Since $(v,s)$ covers only $(v,1)$ and $(vs,s)$ we can also think of $\eta$ as deleting $(v,s)$ and identifying 
        $(v,1)$ with $(vs,s)$.

\begin{figure}[ht]
\caption[The map $\eta:\lbrack 1,srt\rbrack \times \lbrack 1,s\rbrack \rightarrow \lbrack 1,srts\rbrack$.]{The map $\eta:[1,srt] \times [1,s] \rightarrow [1,srts]$, where $[1,srt]$ and $[1,srts]$ are intervals in $(S_4,\set{r,s,t})$.
All elements $(u,v)$ map to $uv$ except $(s,s)$, which maps to $s$.}
\label{fig3}

\centerline{\epsfbox{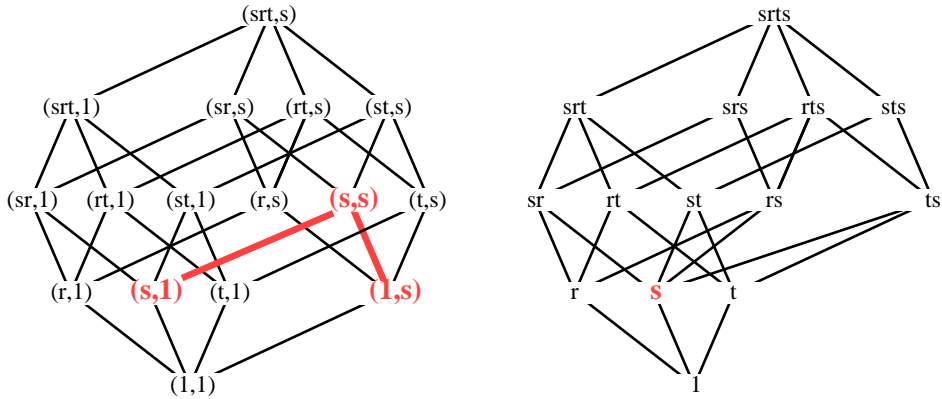}}   
\end{figure}

\begin{figure}[ht]
\caption{The CW-spheres associated to the posets of Figure~\ref{fig3}.}
\label{fig4}

\centerline{\epsfbox{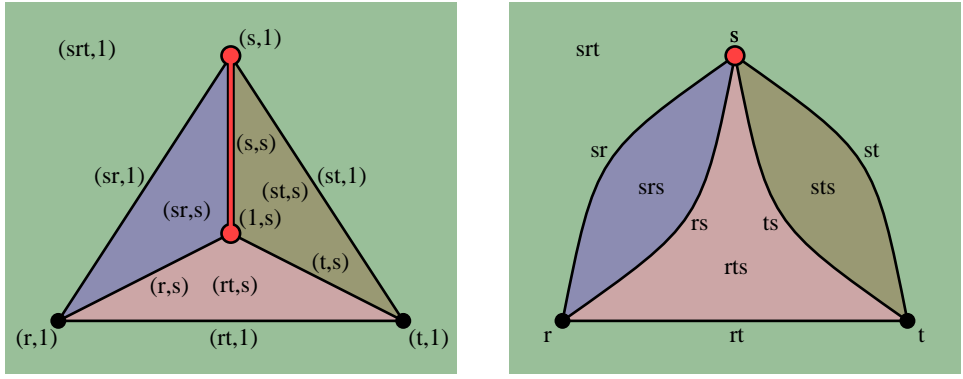}}   
\end{figure}

The map $\eta$ induces a map (also called $\eta$) on the CW-spheres associated to Bruhat intervals, as illustrated in Figure \ref{fig4}.

\begin{prop}
\label{usNotInInterval}
Let $u <us$, $w < ws$ and $us \not\le w$.
Then $vs>v$ for all $v\in [u,v]$, and $\eta$ is an isomorphism.
\end{prop}

\begin{proof}
Suppose for the sake of contradiction that there is a $v\in [u,w]$ with $vs<v$.
Since $us>u$ and $u\le v$, by lifting, $v\ge us$.
By transitivity, $w\ge us$. 
This contradiction shows that that $vs>v$ for all $v\in [u,v]$.

Now, looking back at the proof of Proposition \ref{eta}, we see that
        \[\eta^{-1}(v) = \left\{\begin{array}{ll}
                        \set{(v,1)}   & \mbox{if } v < vs\\
                        \set{(vs,s)}  & \mbox{if } v > vs, \end{array} \right.\]
because in the $v > vs$ case, the other two possible elements of $\eta^{-1}(v)$ don't exist.
Thus the map $\nu$ is an order-isomorphism and therefore $\eta = \bar{\eta}\circ \nu$ is also an order-isomorphism.
\end{proof}

The following corollary is easy.
\begin{corollary}
\label{zeta}
If $u <us$, $w < ws$ and $us \not\le w$ the map $\zeta:[u,w]\rightarrow [us,ws]$ with $\zeta(v) = vs$ is an isomorphism.
\qed
\end{corollary}

\begin{figure}[ht]
\caption[The construction of $\Sh_s(\lbrack 1,rst\rbrack)$ from $\lbrack 1,rst \rbrack \times\lbrack 1,s\rbrack$]{The construction of $\Sh_s([1,rst])$ from $[1,rst]\times[1,s]$, 
where $[1,srt]$ is an interval in $(S_4,\set{r,s,t})$.
The posets are $[1,rst]$, $[1,rst]\times[1,s]$, the same poset with $(s,1),(1,s),(rs,s)$ zipped, and $\Sh_s([1,rst])$.}
\label{fig5}

\centerline{\epsfbox{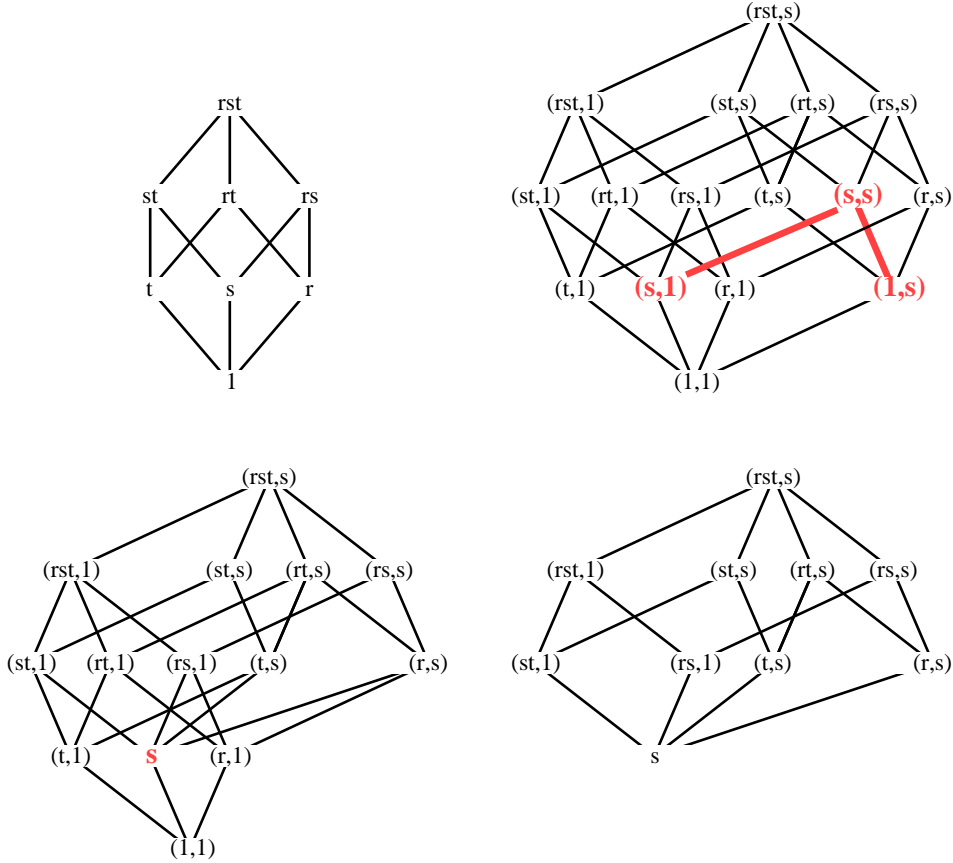}}   

\end{figure}

\begin{figure}[ht]
\caption{The CW-spheres associated to the posets of Figure \ref{fig5}.}
\label{fig6}

\centerline{\epsfbox{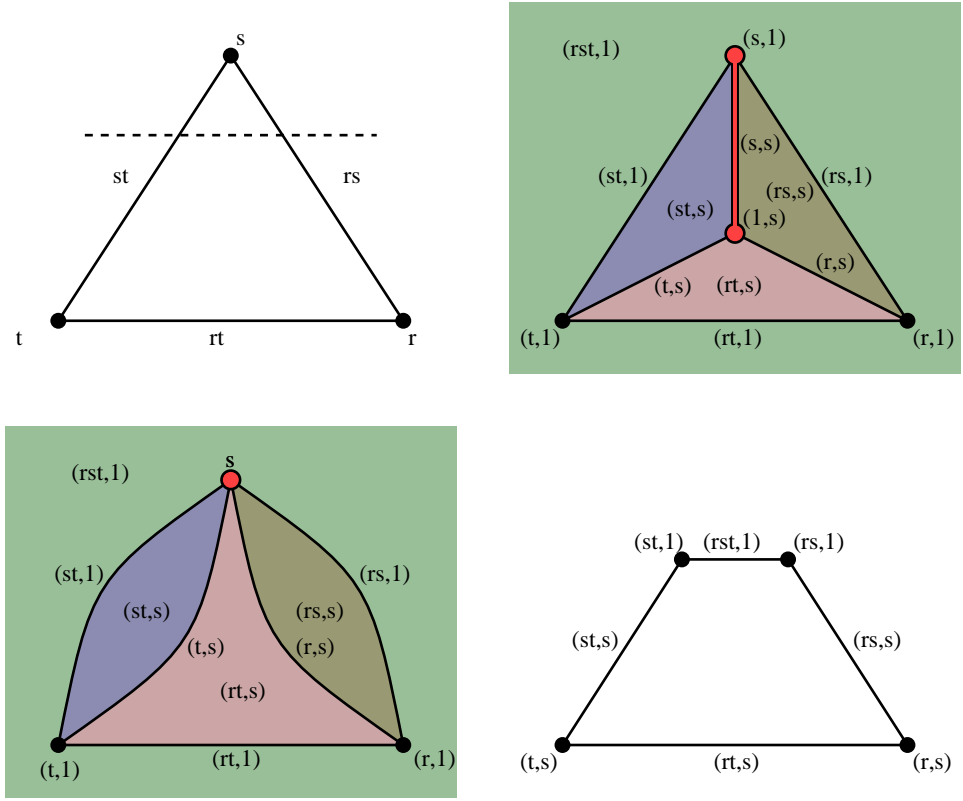}}   

\end{figure}

We would also like to relate the interval $[us,ws]$ to $[u,w]$ in the case where $us \le w$.
To do this, we need an operator on posets corresponding to vertex-shaving on polytopes or CW spheres.
Let $P$ be a poset with $\hat{0}$ and $\hat{1}$, and let $a$ be an atom of $P$.
The {\em shaving} of $P$ at $a$ is an induced subposet of $P \times [\hat{0},a]$ given by:
        \[\Sh_a(P) = \left(\left(P-\set{\hat{0},a}\right)\times \set{a} \right) \cup 
                \left((a,\hat{1}] \times \set{\hat{0}} \right) \cup \set{(\hat{0},\hat{0})}.\]
We can also describe $\Sh_v(P)$ as follows:
Let $P'$ be obtained from $P\times[\hat{0},a]$ by zipping the zipper $((a,\hat{0}),(\hat{0},a),(a,a))$.
Denote by $a$ the element created by the zipping.
Then $\Sh_v(P)$ is the interval $[a,(\hat{1},a)]$ in $P'$.
Figures \ref{fig5} and \ref{fig6} illustrate the operation of shaving.

Let $s \in S$, $u<us$, $w<ws$ and $us \le w$.
Define a map $\theta:\Sh_{us}[u,w]\rightarrow [us,ws]$ as follows.
Starting with $[u,w] \times [1,s]$, zip $((us,1),(u,s),(us,s))$, call the new element $us$, and identify $S_{us}[u,w]$ with the interval $[us,(w,s)]$ in the 
zipped poset.
Now define:
        \begin{eqnarray*}
                \theta(us) & = & us\\
                \theta(v,1) & = & v \mbox{ if } v \in (us,w]\\
                \theta(v,s) & = & \left\{\begin{array}{ll}
                                        vs & \mbox{if  } vs > v\\
                                        v  & \mbox{if  } vs < v. \end{array} \right. \end{eqnarray*}
To check that $\theta$ is well-defined, begin by noting that $\theta(us) \in [us,ws]$, and if $v \in (us,w]$, then $\theta(v,1) = v \in [us,ws]$.
If $v \in [u,w] - \set{u,us}$, there are two possibilities, $\theta(v,s)=vs>v$ or $\theta(v,s)=v>vs$.
In either case, $us<\theta(vs)<ws$ by lifting.
So $\theta$ is well-defined.

\begin{prop}
\label{theta}
The map $\theta:\Sh_{us}[u,w]\rightarrow [us,ws]$ is an order-projection.
\end{prop}
\begin{proof}
Notice that $\theta$, restricted to $\Sh_{us}[u,w] - \set{us}$ is just $\eta$ restricted to an induced subposet.
Recall that in the proof of Proposition \ref{eta}, it was shown that for $v\in [u,ws]$, if $v<vs$ then $(v,1)\in \eta^{-1}(v)$ and if $v>vs$ then 
$(vs,s)\in\eta^{-1}(v)$.
The existence of these elements of $\eta^{-1}(v)$ was used to check that $\eta$ is an order-projection.
The same argument accomplishes most of the present proof.
For $us<x_1\le x_2\le w$, we are done, and it remains to check that for $us\le x\le w$, there exist elements $a\le b$ in $\Sh_{us}[u,w]$ with $\theta(a)=us$ and 
$\theta(b)=x$.
This is easily accomplished by setting $a=us\in\Sh_{us}[u,w]$ and letting $b$ be an element of $\eta^{-1}(x)=\theta^{-1}(x)$. 
If $x=us$, then set $b=us\in\Sh_{us}[u,w]$.
\end{proof}
Figures \ref{fig7} and \ref{fig8} illustrate the map $\theta$ and the corresponding map on CW spheres.

\begin{figure}[ht]
\caption[The map $\theta:\Sh_s(\lbrack 1,rst\rbrack) \rightarrow \lbrack s,rsts\rbrack$.]{The map $\theta:\Sh_s([1,rst]) \rightarrow [s,rsts]$, where $[1,srt]$ and $[s,rsts]$ are intervals in $(S_4,\set{r,s,t})$.
All elements $(u,v)$ map to $uv$ except $(rs,s)$, which maps to $rs$.}
\label{fig7}

\centerline{\epsfbox{bruhat_fig7.ps}}   

\end{figure}

\begin{figure}[ht]
\caption{The CW-spheres associated to the posets of Figure~\ref{fig7}.}
\label{fig8}

\centerline{\epsfbox{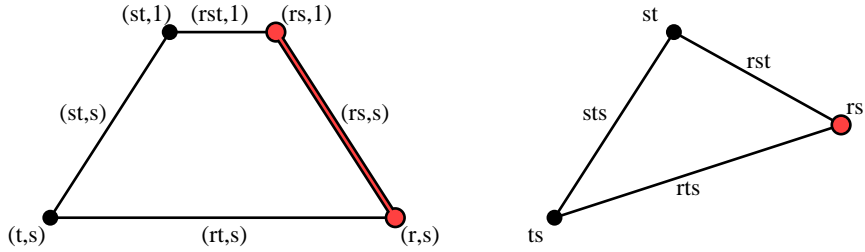}}   

\end{figure}

Since $\eta$ is an order-projection, $[u,ws]$ is isomorphic to the fiber poset of $[u,w]\times[1,s]$ 
        with respect to $\eta$.
Similarly, $[us,ws]$ is isomorphic to the fiber poset of $\Sh_{us}[u,w]$ with respect to $\theta$.
We will show that one can pass to the fiber poset in both cases by a sequence of zippings.

Order the set $\set{v\in (u,w): vs< v}$ linearly such that the elements of rank $i$ in $[u,w]$ precede the elements of rank $i+1$ for all $i$.
Write this order as $v_1,v_2,\ldots,v_k$. 
Define $P_0 = [u,w]\times[1,s]$ and inductively define $P_i$ to be the poset obtained by zipping $((v_i,1), (v_is,s), (v_i,s))$ in $P_{i-1}$.
We show inductively that this is indeed a proper zipping.
First, notice that $(v_i,1)$, $(v_is,s)$ and $(v_i,s)$ are indeed elements of $P_{i-1}$.
The element $(v_i,s)$ has not been deleted yet, and we have not identified $(v_i,s)$ with any element because it is at a rank higher than we have 
yet made identifications.
The only elements ever deleted are of the form $(x,s)$ where $x>xs$, so $(v_i,1)$ and $(v_is,s)$ have not been deleted.
The only identification one could make involving $(v_i,1)$ and $(v_is,s)$ is to identify them to each other,
        and that has not happened yet.
We check the properties in the definition of a zipper:
By Corollary~\ref{thin} and induction, $P_{i-1}$ is thin, so by Proposition~\ref{simplezipper} it is enough to check that Properties (i) and (ii) hold.
Properties (i) and (ii) hold in $P_0$ and therefore in $P_{i-1}$ because we have made no 
        identifications involving $(v_i,1)$, $(v_is,s)$ and $(v_i,s)$ or higher-ranked elements.
The zipper is proper because $(v_i,s)<(w,s)$ in $P_0$ and thus in $P_{i-1}$.

By definition, $\Sh_{us}[u,w]$ is an interval in the poset $P_1$ defined above.
Specifically, $P_1$ was obtained from $[u,w]\times [1,s]$ by zipping $((us,1), (u,s), (us,s))$.
Let $us$ be the element of $P_1$ resulting from identifying $(us,1)$ with $(u,s)$.
Then $\Sh_{us}[u,w]$ is isomorphic to the interval $[us,(w,s)]$ in $P_1$.
The remaining deletions and identifications in the map $\theta$ are really zippings in the $P_i$.  
Therefore they are zippings in the $P_i$ restricted to $[us,(w,s)]$.

We have proven the following:

\begin{theorem}
\label{etathetazippers}
Let $ws>w$, $us>u$ and $u\le w$.
If $us\not\in[u,w]$ then $[u,ws]\cong[u,w]\times[1,s]$ and $[us,ws]\cong[u,w]$.
If $us\in[u,w]$, then $[u,ws]$ can be obtained from $[u,w]\times[1,s]$ by a sequence of zippings, and $[us,ws]$ can be obtained from 
$\Sh_{us}[u,w]$ by a sequence of zippings.
\qed
\end{theorem}

\begin{cor}
\label{spheres}
Bruhat intervals are PL spheres.
\end{cor}

\begin{proof}
One only needs to prove the corollary for lower intervals, because by Proposition \ref{interval} it will then hold for all intervals. 
Intervals of rank 1 are empty spheres.
It is easy to check that a lower interval under an element of rank 2 is a PL 0-sphere.
Given the interval $[1,w]$ with $l(w)\ge 3$, there exists $s\in S$ such that $ws<w$.
Then $[1,w]$ can be obtained from $[1,ws]\times [1,s]$ by a sequence of zippings.
By induction $[1,ws]$ is a PL sphere, and thus $[1,ws]\times [1,s]$ is as well.
By repeated applications of Theorem \ref{pl}, $[1,w]$ is a PL sphere.
\end{proof}

The following observation will be helpful in Section \ref{recursion section}, when Theorem \ref{etathetazippers} is combined with Theorem \ref{cd thm}.
\begin{prop}
\label{intervals}
For $1\le i \le k$,
        \[ [(u,1),(v_i,1)]_{P_{i-1}} \cong [(u,1),(v_i,1)]_{P_0}\]
        \[ [(v_i,s),(w,s)]_{P_{i-1}} \cong [(v_i,s),(w,s)]_{P_0}\]
\end{prop}

\begin{proof}
The second statement is obvious because of the way the $v_j$ were ordered.
For the first statement, there is the obvious order-preserving bijection between the two intervals.
The only question is whether the right-side has any extra order relations.  
Extra order relations will occur if for some $v$ with $v>vs$, there exists $(x,1)$ with
        $(x,1)\le (vs,s)$ but $(x,1)\not\le (v,1)$.
This is ruled out by transitivity.
\end{proof}

\section[cd-Index recursion]{A Recursion for the cd-index of Bruhat intervals}
\label{recursion section}
Theorems \ref{cd thm} and \ref{etathetazippers} yield Theorem \ref{formula}, a set of recursions for the cd-indices of Bruhat intervals.
In this section we prove Theorem \ref{formula}, then apply it to determine the affine span of the cd-indices of Bruhat intervals.

For $v\in W$ and $s\in S$, define $\sigma_s(v) := l(vs)-l(v).$
Thus $\sigma_s(v)$ is 1 if $ v$ is lengthened by $s$ on the right and $-1$ if $v$ is shortened by $s$ on the
        right. 

\begin{theorem}
\label{formula}
Let $u<us$, $w<ws$ and $u\le w$.

If $us\not\in[u,w]$, then $\Psi_{[u,ws]}=\Pyr\Psi_{[u,w]}$, and $\Psi_{[us,ws]}=\Psi_{[u,w]}$.

If $us\in[u,w]$, then
\begin{eqnarray*}
\Psi_{[u,ws]} &=& \Pyr\Psi_{[u,w]} - \sum_{v\in(u,w):vs< v} \Psi_{[u,v]}\cdot d\cdot\Psi_{[v,w]} \\
              &=& \frac{1}{2}\left(\Psi_{[u,w]}\cdot c + c\cdot\Psi_{[u,w]} + \sum_{v\in (u,w)}
                        \sigma_s(v)\Psi_{[u,v]}\cdot d \cdot \Psi_{[v,w]} \right)\\
\Psi_{[us,ws]} &=& \Sh_{us}\Psi_{[u,w]} - \sum_{v\in (us,w):vs< v} \Psi_{[us,v]}\cdot d\cdot\Psi_{[v,w]} \\
               &=& \Psi_{[u,w]} + \frac{1}{2}\left(\Psi_{[us,w]}\cdot c - c\cdot\Psi_{[us,w]} + 
                        \!\!\sum_{v\in (us,w)} \sigma_s(v) \Psi_{[us,v]}\cdot d \cdot \Psi_{[v,w]} \right)
\end{eqnarray*}
\end{theorem}

The first line of each formula looks like an augmented coproduct \cite{Ehr-Fox} on a Bruhat interval, with an added sign.
The second line of each formula is more efficient for computation, because the formulas in Proposition \ref{cd derivation} 
are more efficient than the forms quoted in Proposition \ref{cd ops}.

\begin{proof}[Proof of Theorem \ref{formula}]
The statement for $us\not\in[u,w]$ follows immediately from Propositions \ref{usNotInInterval} and \ref{zeta}.
Define the $P_i$ as in Section \ref{Building}.
Thus by Theorem \ref{cd thm}, 
 \[\Psi_{P_{i-1}}-\Psi_{P_{i}} = \Psi_{[(u,1),(v_i,1)]_{P_{i-1}}}\cdot d\cdot\Psi_{[(v_i,s),(w,s)]_{P_{i-1}}}.\]
Since $P_k = [u,ws]$, sum from $i=1$ to $i=k$ to obtain
    \[\Psi_{[u,ws]}=\Psi_{P_0}-\sum_{j=1}^k \Psi_{[(u,1),(v_j,1)]_{P_{j-1}}}\cdot d\cdot\Psi_{[(v_j,s),(w,s)]_{P_{j-1}}}.\]
By Proposition \ref{intervals}, $[(u,1),(v_j,1)]_{P_{j-1}} \cong [(u,1),(v_j,1)]_{P_0}$.
This in turn is isomorphic to $[u,v_j]$.
Similarly, by Proposition \ref{intervals}, $[(v_j,s),(w,s)]_{P_{j-1}} \cong [(v_j,s),(w,s)]_{P_0}$ which is isomorphic to $[v_j,w]$.
Thus we have established the first line of the first formula in Theorem \ref{formula}.
The second line follows from the first by Proposition \ref{cd ops}.

A similar proof goes through for the second formula.
The isomorphisms from Proposition \ref{intervals} restrict to the appropriate isomorphisms for the $[us,ws]$ case.
The second line of the formula follows by Proposition \ref{cd ops}.
\end{proof}

In \cite{Bayer-Billera}, Bayer and Billera show that the affine span of the cd-indices of polytopes is the entire affine space of monic 
        cd-polynomials.
As an application of Theorem \ref{formula}, we prove that the cd-indices of Bruhat intervals have the same affine span.

\begin{theorem}
\label{span}
The set of cd-indices of Bruhat intervals spans the affine space of cd-polynomials.
\end{theorem}
\begin{proof}
The space of cd-polynomials of degree $n-1$ has dimension $F_n$, the Fibonacci number, with $F_1=F_2=1$ and 
        $F_n=F_{n-1}+F_{n-2}$.
For each $n$ we will produce a set $\F_n$ consisting of $F_n$ reduced words, corresponding to group elements whose lower Bruhat intervals 
have linearly independent cd-indices.

Let $(W,S:=\set{s_1,s_2,\ldots})$ have a complete Coxeter graph with each edge labeled~3.
Each $\F_n$ is a set of reduced words of length $n$ in $W$, with $\F_1=\set{s_1}$, $\F_2=\set{s_1s_2}$ and 
        \[\F_n=\F_{n-1}s_n \,\cup\!\!\!\!\cdot \,\,\,\,s_n\F_{n-2}s_n,\]
where $\cup\!\!\!\cdot \,\,\,$ means disjoint union.
Given a word $w\in\F_{n-1}$, by Proposition \ref{usNotInInterval}, $[1,ws_n]\cong\Pyr[1,w]$, so 
        $\Psi_{[1,ws_n]}=\Pyr(\Psi_{[1,w]})$.
Similarly, given a word $w'\in\F_{n-2}$, $[1,s_nw]\cong\Pyr[1,w]$.
Since $s_n$ does not commute with any other generator, and since $s_n$ is not a letter in $w$, by Proposition \ref{eta}, $[1,s_nws_n]$ is obtained from 
$[1,s_nw]$ by a single zipping.
In particular, by Theorem \ref{formula}, $\Psi_{[1,s_nws_n]}=\Pyr^2(\Psi_{[1,w]})-d\cdot\Psi_{[s_n,s_nw]}$.
By a ``left'' version of Corollary \ref{zeta}, we have $[s_n,s_nw]\cong[1,w]$, so $\Psi_{[1,s_nws_n]}=\Pyr^2(\Psi_{[1,w]})-d\cdot\Psi_{[1,w]}$.
Let $\Psi(\F_n)$ be the set of cd-indices of lower intervals under words in $\F_n$.
The proof is now completed via Proposition~\ref{indep}, below, which in turn depends on Lemma~\ref{Pyr}.
\end{proof}

\begin{prop}
\label{indep}
For each $n\ge 1$, the $F_n$ cd-polynomials in $\Psi(\F_n)$ are linearly independent.
\end{prop}

\begin{proof}
As a base for induction, the statement is trivial for $n=1,2$.
For general $n$, form the matrix $M$ whose rows are the vectors $\Psi(\F_n)$, written in the cd-index basis.
Order the columns by the lexicographic order on cd-monomials.
Order the rows so that the cd-indices in $\Psi(\F_{n-1}s_n)$ appear first.
We will show that there are row operations which convert $M$ to an upper-unitriangular matrix.
Notice that for each $w\in \F_{n-2}$, $\Pyr^2\Psi_{[1,w]}$ occurs in $\Psi(\F_{n-1}s_n)$.
Also, $\Pyr^2\Psi_{[1,w]}-d\cdot\Psi_{[1,w]}$ occurs in $\Psi(s_n\F_{n-2}s_n)$.
Thus by row operations one obtains a matrix $M'$ whose rows are first $\Pyr(\Psi(\F_{n-1}))$, then $d\cdot\Psi(\F_{n-2})$.
By induction, there are row operations which convert the matrix with rows $\Psi(\F_{n-1})$ to an upper-unitriangular matrix.
By Corollary~\ref{Pyr}, these yield row operations which give the first $F_{n-1}$ rows of $M'$ an upper-unitriangular form.
Also by induction, there are row operations which convert the matrix with rows $\Psi(\F_{n-2})$ to an upper-unitriangular matrix.
Corresponding operations applied to the rows $d\cdot\Psi(\F_{n-2})$ of $M'$ complete the reduction of $M'$ to upper-unitriangular form.
\end{proof}

The proof of Theorem~\ref{span} uses infinite Coxeter groups.
It would be interesting to know whether the cd-indices of Bruhat intervals in finite Coxeter groups also span, and whether a spanning set of 
intervals could be found in the finite Coxeter groups of type~A.

\section[Bounds]{Bounds on the cd-index of Bruhat intervals}
\label{bounds section}
In this section we discuss lower and upper bounds on the coefficients of the cd-index of a Bruhat interval.
The conjectured lower bound is a special case of a conjecture of Stanley \cite{Stanley2}.
\begin{conj}
\label{nonneg}
For any $u \le v$ in $W$, the coefficients of $\Phi_{[u,v]}$ are non-negative.
\end{conj}
The coefficient of $c^n$ is always 1, and for the other coefficients the bound is sharp because the dihedral group $I_2(m)$ has cd-index $c^m$. 
Computer studies have confirmed the conjecture in $S_n$ with $n\le 6$.

The conjectured upper bounds are attained on Bruhat intervals which are isomorphic to the face lattices of convex polytopes.
For convenience, we will say that such intervals ``are'' polytopes.
We now use the results of Section~\ref{Building} to construct these intervals.

A polytope is said to be {\em dual stacked} if it can be obtained from a simplex by a series of vertex-shavings.
As the name would indicate, these polytopes are dual to the {\em stacked polytopes} \cite{Klei-Lee} which we will not define here. 
There are Bruhat intervals which are dual stacked polytopes.
Let $W$ be a universal Coxeter group with Coxeter generators $S:=\set{s_1,s_2,\ldots s_{d+1}}$.
Define $C_k$ (for ``cyclic word'') to be the word $s_1s_2\cdots s_k$, where the subscript $k$ is understood to mean $k~(\mbox{mod }~d+1)$.
So for example, if $d=2$, then $C_7=s_1s_2s_3s_1s_2s_3s_1$.
In a universal Coxeter group, every group element corresponds to a unique reduced word. 
Thus we will use these words interchangeably with group elements.

\begin{prop}
\label{dual stacked}
The interval $[C_k,C_{d+k+1}]$ in $W$ is a dual stacked polytope of dimension $d$ with $d+k+1$ facets.
\end{prop}
\begin{proof}
The proof is by induction on $k$.
For $k=0$, the interval is $[1,s_1s_2\cdots s_{d+1}]$, and because all subwords of $s_1s_2\cdots s_{d+1}$ are distinct group elements, this interval is a 
Boolean algebra---the 
face poset of a $d$ dimensional simplex.
For $k>0$, the interval $[C_k,C_{d+k+1}]$ is obtained from the interval $[C_{k-1},C_{d+k}]$ by shaving the vertex $C_{k}$ and then possibly by performing a sequence of zippings.
By induction, $[C_{k-1},C_{d+k}]$ is a dual stacked polytope of dimension $d$ with $d+k$ facets.
The zippings correspond to elements of $(C_{k},C_{d+k})$ which are shortened on the right by $s_{d+k+1}$.
We will show that in fact there are no such elements.
Let $v\in(C_{k},C_{d+k})$, and let $v$ also stand for the unique reduced word for the element $v$.
Since every element of $W$ has a unique reduced word, the fact that $C_k< v$ means that $v$ contains $C_k$ as a subword.
But the only subword $C_k$ of $C_{d+k}$ is the first $k$ letters of $C_{d+k}$.
Thus $v$ is a subword of $C_{d+k}$ consisting of the first $k$ letters and at least one other letter.
Now, $s_k\not\in\set{s_{k+1},s_{k+2},\ldots,s_{d+k}}$, so $v$ ends in some generator other than $s_k$, and therefore, $v$ is not shortened on the right by 
$s_k=s_{d+k+1}$.
Thus there are no zippings following the shaving, and so $[C_k,C_{d+k+1}]$ is a dual stacked polytope of dimension $d$ with $d+k+1$ facets.
\end{proof}

Proposition~\ref{dual stacked} shows that the following conjectured upper bound is sharp if indeed it holds.

\begin{conj}
\label{dual stacked upper}
The coefficientwise maximum of all cd-indices $\Phi_{[u,v]}$ with $l(u)=k$ and $l(v)=d+k+1$ is attained on a Bruhat interval which is 
isomorphic to a dual stacked polytope of dimension $d$ with $d+k+1$ facets. 
\end{conj}

This conjecture is natural in light of Proposition \ref{dual stacked} and Theorem \ref{formula}.
There are two issues which complicate the conjecture.
First, any proof using Theorem \ref{formula} requires non-negativity (Conjecture \ref{nonneg}).
Second, and perhaps even more serious, there is the issue of commutation of operators.

Given $p\in P$ denote the corresponding ``downstairs'' element in $\Pyr(P)$ by $p$ and the ``upstairs'' element by $p'$.
Denote the operation of zipping a zipper $(x,y,z)$  by ${\Z}_z$.
Then $\Pyr{\Z}_zP\cong {\Z}_{z'}{\Z}_z\Pyr P$.
The triple $(x',y',z')$ becomes a zipper only after ${\Z}_z$ is applied.
Pyramid and shaving also commute reasonably well: $\Pyr\Sh_aP\cong {\Z}_{a'}\Sh_a\Pyr P$.

However, zipping does not in general commute nicely with the operation of shaving off a vertex $a$.
Given $p\ne a\in P$ denote the corresponding element of $\Sh_aP$ again by $p$, and if in addition $p>a$, write $\bar{p}$ for the new 
element created by shaving.
If $z>a$ and $a\not\in\set{x,y}$ then $\Sh_a{\Z}_zP\cong {\Z}_z{\Z}_{\bar{z}} \Sh_aP$.
If $z\not>a$ then $\Sh_a{\Z}_zP\cong {\Z}_z\Sh_aP$.
However, if $x$ and $y$ are vertices then $\Sh_{xy}{\Z}_zP=\Sh_zP$, where $\Sh_z$ is the operation of shaving off the edge $z$.

Since the pyramid operation commutes nicely with zipping, and any lower interval is obtained by pyramid and zipping operations, it is possible 
to obtain any lower interval by a series of pyramid operations followed by a series of zippings.
Thus by Theorem \ref{formula}: 
\begin{theorem}
\label{lower}
Assuming Conjecture \ref{nonneg}, for any $w$ in an arbitrary Coxeter group,
        \[\Phi_{[1,w]} \le \Phi_{B_{l(w)}}.\]
\end{theorem}
Here $B_n$ is the Boolean algebra of rank $n$.
It is not true that the cd-index of general intervals is less than that of the Boolean algebra of appropriate rank.
For example, $[1324,3412]$ is the face lattice of a square, with $\Phi_{[1324,3412]} = c^2+2d$.
However, $\Phi_{B_3} = c^2 + d$.

Equation (\ref{ab formula}) in the proof of Theorem \ref{cd thm} is a formula for the change in the ab-index under zipping.
Thus Theorem \ref{formula} has a flag h-vector version, and since the flag h-vectors of Bruhat intervals are known to be nonnegative, the 
following theorem holds.
\begin{theorem}
\label{boolean h}
For any $w$ in an arbitrary Coxeter group,
        \[\Psi_{[1,w]} \le \Psi_{B_{l(w)}}.\]
\end{theorem}
Here ``$\,\le$'' means is coefficientwise comparison of the ab-indices, or in other words, comparison of flag h-vectors.

\section{Acknowledgments}
The author wishes to thank his adviser, Vic Reiner, for many helpful conversations, as well as Lou Billera, Francesco Brenti, Kyle Calderhead, 
Richard Ehrenborg, Margaret Readdy and Michelle Wachs for helpful comments.

\newcommand{\journalname}[1]{\textrm{#1}}
\newcommand{\booktitle}[1]{\textrm{#1}}

\bibliographystyle{plain}

\end{document}